\documentclass[a4paper,10pt]{article}
\usepackage{amssymb,rotating,graphics,epsfig,float,graphics}
\usepackage{amsmath,verbatim,a4wide}
\usepackage{hyperref}
\usepackage{theorem}
\theoremstyle{break} 
{\theorembodyfont{\rmfamily} }
{\theorembodyfont{\rmfamily} \newtheorem{prop}{Theorem}[section]}

\usepackage{color}

\definecolor{blue}{rgb}{0,0,1}
\definecolor{red}{rgb}{1,0,0}

\begin{document}

\author{Wilberd van der Kallen and Ferdinand Verhulst\\
Mathematisch Instituut, University of Utrecht\\
PO Box 80.010, 3508 TA Utrecht The Netherlands}
\title{Explorations for alternating FPU-chains with large mass}

\maketitle 
 \begin{abstract} 
We show interaction between high- and low-frequency modes in periodic $\alpha$-FPU chains with alternating 
 large masses. This high-low frequency interaction is known if the number of particles $N$ is a 4-fold, 
 our treatment discusses the difficult case where the number of particles $N=2p$ involves $p$ prime. 
 A key role is played by identifying symmetric invariant manifolds, thus reducing the dimension of the problems drastically, 
 and a {\sc Mathematica} programme focused on these systems. We show explicitly high-low frequency interaction for systems  
 with $2pn$ particles where $2 \leq p \leq 47$ is prime and $n$ is an arbitrary natural number. In addition we have 
 strong arguments for interactions in arbitrary large chains.  
 \end{abstract} 
 
 MCS2020 classification: 34-04, 34A26, 34C45, 37J99, 70H33, 70K70. \\

 \begin{quote} 
 My ideas and my opinions feel like groping their way, faltering, stumbling and making false steps; 
 and when I have gone as far as is possible for me, I am not satisfied at all: I see more fields after this, 
 vaguely and in a haze that I cannot clear up.\\
 {\em Michel de Montaigne, Essais vol.~I (1580)}
  \end{quote}

\section{Introduction} 
 The present paper is a sequel to \cite{FVa}. 
In large nonlinear chains and in nonlinear wave theory one of the basic questions is whether interaction between 
very different parts of the spectrum are of importance, see for an example \cite{DLS18}. Fermi-Pasta-Ulam (FPU) chains 
play a paradigmatic part 
in discussions of particle interactions, see the surveys \cite{CRZ05}, \cite{G08} and for a nice introduction to FPU-like 
models \cite{B20}. Originally FPU chains were studied to understand fundamental questions (ergodicity) in statistical 
mechanics but the perspective has changed to bifurcations, periodic solutions and waves in nonlinear chains. 
As usual in theoretical physics and dynamical systems theory in such studies symmetry considerations are basic. 
It came as a surprise that for certain FPU-like chains widely separate parts of the spectrum can interact at 
the nonlinear level.
For  periodic FPU-chains with $N$ alternating large masses it was shown in \cite{FVa}, see also \cite{BVapp}, 
that for  $\alpha$-chains with  
$N=2n$ particles, in particular if $n$ is even ($N$ a 4-fold),  we have significant interactions caused by external 
forcing of the low-frequency acoustic modes by a stable 
or unstable high-frequency optical normal mode. In the case of $n$ prime the analysis was restricted to 
the examples $N= 6, 10$; 
this was caused by the formidable problems posed by linear algebra 
manipulations to transform to quasi-harmonic (normal mode) equations. For the present paper
a  special {\sc Mathematica} programme was developed and used to extend the analysis to periodic chains 
with many more particles. \\

It was shown in \cite{FVa} that in the case $N= 2p$ with $p$ odd invariant 
manifolds exist allowing drastic reduction of the number of degrees-of-freedom (dof) to $p-1$. The dynamics is 
more complicated 
than in the case $N$ with $p$ even. The novel results in the present paper are both quantitative and 
qualitative. 
For $p$ odd and $3 \leq p \leq 47$ we can extend the analysis to find significant 
 interactions between widely different parts of the spectrum, but even more interestingly, we find 
 a surprising circular dependence of the equations of motion involving the interactions.
 
In general the
 existence and use of invariant manifolds for FPU-chains is very important because it produces insight in the global 
dynamics of FPU-chains and enables us to replace the problems by studying submanifolds with less degrees-of-freedom (dof). 
In a number of seminal papers
Czechin et al. identified {\em bushes} of solutions forming submanifolds for classical (equal masses) 
FPU-chains; see \cite{C05} and further references there; the paper also generalised the results found for two-mode invariant 
manifolds in \cite{PR}. Our analysis is related to this approach, also to ``localisation analysis'' for FPU-chains in \cite{CEB}.  
See also the analysis of acoustic and optical vibrations for a monatomic lattice in \cite{ADR}.\\

The spatially periodic FPU-chain with $N$ particles where the first oscillator is connected with the last one can be described
by the Hamiltonian 
 \begin{equation} \label{Hfpu}
 H(p, q) = \sum_{j=1}^N \left( \frac{1}{2m_j}p_j^2 + V(q_{j+1} - q_j)\right).
 \end{equation} 

For a general introduction to Hamiltonian dynamics see \cite{BS12}. Following  \cite{GGMV},
we choose the number $N=2n$ of particles even and take the odd masses $m_{2j+1}$ equal to $1$, the much larger 
even masses 
$m_{2j}=m=\frac{1}{a}$, where $a>0$ is small. This chain is an example of an alternating FPU-chain. With this choice of the  
masses, 1 and $m$, the eigenvalues of the  PFU equations of motion near equilibrium split in 2 groups of $n$ eigenvalues
with size resp.\ $O(a)$ (the acoustic group) and size $2+ O(a)$ (the optical group). The cases of $O(1)$ choices of mass $m$ 
was discussed in \cite{BValt}. 
 
We consider the Hamiltonian near stable equilibrium $p=q=0$, and use a potential $V$ of the form 
\[ V(z) = \frac{1}{2} z^2+\frac{\alpha}{3} z^3+\frac{\beta}{4}z^4\,,\]
and speak of an $\alpha$-chain if $\alpha \neq 0$, $\beta=0$ and of a $\beta$-chain if $\alpha =0$, $\beta \neq 0$. 
In this study we will discuss only $\alpha$-chains.

\subsection{Two reduction methods}
 It was shown in \cite{BValt} and 
used in \cite{FVa} 
that an {\em embedding theorem} is useful. 
\begin{prop} \label{embed}
Consider the equations of motion induced by Hamiltonian \eqref{Hfpu} for $\alpha \beta \neq 0$ and 
$\alpha$- or $\beta$-chains, with alternating masses 
$1, m >0$ and $n$ (even) particles. Suppose $k$ is a multiple of $n$ and consider the equations of motion induced by 
Hamiltonian \eqref{Hfpu} with identical $\alpha$, $\beta$, $m$  and $2kn$ dof, then there exists a restriction of this larger 
Hamiltonian system that is equivalent to the first system with $2n$ dof. 
\end{prop} 
So we have that each alternating periodic FPU-chain with $2n \geq 4$ particles 
occurs isomorphically as an invariant submanifold in all subsystems with $2kn$ particles ($k = 2, 3, \ldots$). 
This makes the study of small alternating FPU-chains relevant for larger alternating systems.\\
 If we have studied for instance an alternating chain with 6 particles, we will find the dynamics of this chain in a submanifold 
 of alternating chains with $N= 12$, $18$, $24$ etc.~particles. From a reverse point of view, if we study an alternating chain with 
 32 particles we can at least expect to find invariant manifolds with dynamics of 4, 8 and 16 particles. \\ 
 For $N=4$ high-low frequency interaction was shown in \cite{FVa}, so we will have this interaction in all alternating 
 FPU-chains with $N$ a 4-fold. What remains are the infinite number of cases where $N=2p$ with $p$ odd. 
 If $p$ is odd we can factorise in a product of prime numbers. Showing interaction for one of these prime 
numbers implies interaction for this odd $p$ because of the embedding theorem. So we have to consider only the cases 
of $p$ prime. \\

To understand more of the dynamics in the case $N=2p$ with $p$ prime we have a second reduction method based on
the existence of certain lower-dimensional invariant manifolds derived in \cite{FVa}. 
As we will show in section \ref{sec2} in this case there exist 2 symmetric invariant manifolds corresponding with 
$(p-1)$ dof. Using this reduction  we will study high-low frequency interaction in alternating chains for $4 \leq N \leq 104$. 

\subsection{Set-up of the paper} 
As mentioned above it was shown in \cite{FVa} that with  $N$ a 4-fold and $N=2p$, $p= 3$, $5$, there is 
considerable interaction between the 
acoustic and optical group; the analysis was mainly based on asymptotics (averaging). The embedding theorem \ref{embed} 
implies that this high-low frequency interaction holds in system with $N$ a multiple of these values.  
In sections \ref{sec3}-\ref{sec4} we 
extend this analysis with stability considerations and aspects of normal modes. \\
In section \ref{sec5} we consider $p=9$ (18 particles chain) to obtain insight in the case where $p$ is odd but not prime. 
Interestingly  we find after reduction to symmetric invariant manifolds  
eight 2nd order equations with the following dependence: 
\[ x_1 \rightarrow x_8 \rightarrow x_5 \rightarrow x_2 \rightarrow x_7 \rightarrow x_6 \rightarrow x_1, \] 
meaning that acoustic $x_1^2$ forces optical $x_8$,  $x_8^2$ forces acoustic $x_5$ etc. \\ 
For $N =2p$ with $5\leq p \leq47$ ($p$ prime and odd) we will consider the results of interaction in section \ref{sec6} 
using a {\sc Mathematica} programme without giving technical details. \\
One remarkable feature is that applying the reduction method for the problem with $(p-1)$ dof, we obtain 
as in the case $p=9$ ``forcing squares'' 
as nonlinearities with one of the acoustic group and one optical. This involves a surprising circular 
dependence of the equations of motion on each other as shown explicitly for $p=9$. For the discussion see 
subsection \ref{forcing}, where we also consider forcing of optical by optical or acoustic by acoustic, which 
hides the circularity a bit.
Finally two cartoon problems are discussed to demonstrate the influence of mixed quadratic terms on 
the right-hand sides of the equations of motion.

\medskip 
The figures were produced by using the programme {\sc Matcont} under {\sc Matlab}; the numerics utilises ode78 with (usually) 
absolute and relative accuracy $10^{-10}$. The {\sc Mathematica} notebook performing the transformations from equations 
of motion to quasi-harmonic equations for normal modes is described in section \ref{sec6} and \cite{KVN20}.

\section{Models with $p$ prime} \label{sec2} 
From \cite{FVa} we have that 
we can apply symmetries to the system induced by Hamiltonian \eqref{Hfpu} with $2p$ particles, $p \geq 3$ prime. We assume: 
\begin{eqnarray} \label{IM2ps} 
\begin{cases} 
q_p(t)  = q_{2p}(t)=0, \\
q_2(t)  = -q_{2p-2}(t),~ q_4(t)= -q_{2p-4}(t), \ldots q_{p-1}(t) = -q_{p+1}(t),\\ 
q_1(t)  = -q_{2p-1}(t),~q_3(t)= -q_{2p-3}(t), \ldots q_{p-2}(t) = -q_{p+2}(t). 
\end{cases}
\end{eqnarray} 
The symmetry assumptions imply that the value of the momentum integral (see \cite{FVa}) vanishes. 
The Hamiltonian system with $2p$ dof contains a submanifold of $2p-2$ dof which is a 4-fold;  we are left with two $(p-1)$ dof 
systems with identical dynamics, quite a reduction. The $(p-1)$ dof system is of the form: 
\[ \ddot{q} + B q = \alpha N(q) \] 
with $q= (q_1, \ldots, q_{p-1})^T,\, B$ is a $(p-1) \times (p-1)$ matrix and $N(q)$ a homogeneous vector, quadratic in the 
$q$ variables. Explicitly for $p > 5$ :
\begin{eqnarray} \label{IM2p} 
 \begin{cases}
 \ddot{q}_1 +2 q_1  -q_2 & = \alpha(q_2^2 -2q_1q_2),\\
m \ddot{q}_2 +2 q_2 - q_{1} -q_3 & = \alpha(q_3^2 -2q_2q_3+2q_1q_2 -q_{1}^2),\\
 \ddot{q}_3 +2 q_3 - q_{2} -q_4 & = \alpha(q_4^2 -2q_3q_4 +2q_2q_3 -q_{2}^2),\\ 
m \ddot{q}_4 +2 q_4 - q_{3} -q_5 & = \alpha(q_5^2 - 2q_4q_5 + 2q_3q_4-q_{3}^2),\\ 
 \cdots & = \alpha (\cdots ),\\ 
m \ddot{q}_{p-1} +2 q_{p-1} - q_{p-2}  & = \alpha(2q_{p-1}q_{p-2} -q_{p-2}^2). 
 \end{cases}
\end{eqnarray} 
The normal mode frequencies of system \eqref{IM2p} are derived from the single eigenvalues of $B$, resp. 
$\surd{2} + O(a)$ and $O(\surd{a})$.  
To put the system in quasi-harmonic form we use a linear transformation that diagonalises B and puts the eigenvalues 
on the diagonal. The linear transformation of $q$ keeps the nonlinearities quadratic, 
resulting in interactions.  In \cite{FVa} the cases $p= 2, 3, 5$ were discussed and consequently chains with $4n, 6n, 10n$ 
particles, $n$ a natural number.

\section{Periodic FPU $\alpha$-chain with 6  particles} \label{sec3}   

The invariant manifolds described by \eqref{IM2ps} can be used to consider systems with $6$ particles ($p=3$).  
Interactions of this chain were studied in \cite{FVa} using averaging, we add an analysis of the symmetric 
invariant manifolds,  
their stability  in 12-dimensional phase-space and some new aspects of the general dynamics.  The stability analysis of invariant 
manifolds is not easy in the case of high dimensions; we will use the stability results of individual modes and 
in addition we obtain insight 
from application of the Poincar\'e recurrence theorem.
The eigenvalues, producing squared frequencies of the linearised system are: 
\begin{equation}
\omega_1^2= 2(1+a),\ \omega_{2,3}^2= a+1+ \sqrt{a^2-a+1},\ \omega_{4,5}^2= a+1- \sqrt{a^2-a+1},\ \omega_6^2=0. 
\end{equation} 
For $a=0.01$ the corresponding eigenvalues are $2.02$, $2.00504$, $0.0149623$, $0$. 

\subsection{Stability of symmetric invariant manifolds} 

\begin{figure}[ht]
\begin{center}
\includegraphics[width=7cm]{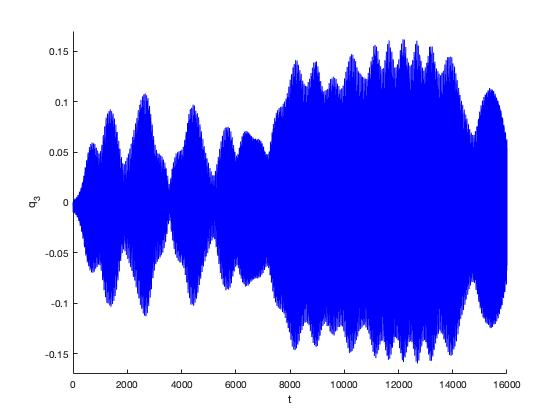}\, \includegraphics[width=7cm]{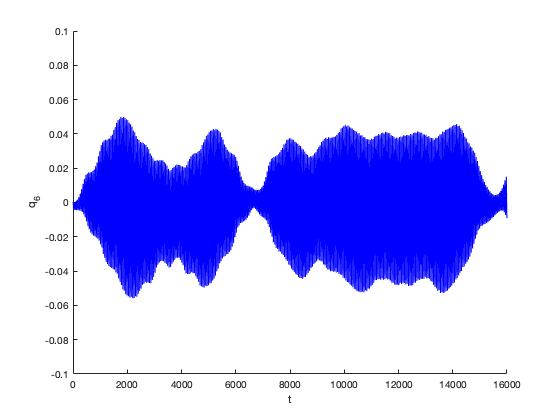}
\end{center}
 \caption{Solutions in the  case of $N=6$ ($p=3$) described by the chain induced by Hamiltonian \eqref{Hfpu}
 close to the symmetric invariant manifolds for 16000 timesteps; 
 $a=0.01$, $\alpha = 1$. The initial values close to symmetry are: $q_1(0)= 0.08$, $q_2(0)= -0.085$, $q_3(0)=0$, 
 $q_4(0)= 0.075$, $q_5(0)= -0.07$, $q_6(0)=0$ and initial velocities zero. 
 We observe forcing of the modes $q_3$, $q_6$ showing instability of the symmetric invariant manifolds. }
\label{figp=3unst}
\end{figure}   

The equations of motion induced by Hamiltonian \eqref{Hfpu} contain symmetric invariant manifolds described by system 
\eqref{IM2p}. In the general  case with the number of particles $N=2p$ with $p$ odd, we have that $q_p(t)=q_{2p}(t) =0$, $t \geq 0$. In the case $p=3$ we have: 
\begin{equation} \label{IM6}
q_3(t)=q_6(t)=0, \, q_2(t)=-q_4(t), \,q_1(t)= -q_5(t). 
\end{equation} 

 It is clear from system \eqref{IM2p} that on breaking the symmetry conditions, in the case $p=3$, 
the modes $q_3$, $q_6$ (in general $q_p$, $q_{2p}$) are forced. So the solutions $q_3=q_6=0$ will not persist, 
the symmetric invariant manifolds are expected to be unstable. Instead of explicit matrix calculations we can demonstrate 
the instability by high precision numerical analysis as indicated in the Introduction. 
We show this explicitly for the case of 6 particles in fig.~\ref{figp=3unst}. \\ 
A related criterion is using the Poincar\'e recurrence theorem. 
Hamiltonian flow on a bounded energy manifold will after some time return arbitrary close to the initial value. Of course the time for this 
return will depend on our definition of closeness and the particular system studied, see for discussion and references \cite{FV18}. 
The computation illustrated in fig.~\ref{figp=3unst} should lead after some time to $q_3(t)$, $q_6(t)$ simultaneously 
approaching zero arbitrarily close but this takes longer than 16000 timesteps. The dynamics shows more complexity than expected. 

\subsection{Symmetric invariant manifold}  
If $2p=6$ the dynamics in the symmetric invariant manifold for $q_1, q_2$ is described by:  
\begin{eqnarray} \label{IMsym6} 
\begin{cases}
\ddot{q}_1 +2q_1 -q_2 & =  \alpha(q_2^2 -2q_1q_2),\\
m \ddot{q}_2+2q_2 -q_1 & =  \alpha(2q_1q_2 - q_1^2). 
\end{cases} 
\end{eqnarray} 
System \eqref{IMsym6} has 4 critical points (we take $\alpha=1$): $(q_1, q_2, \dot{q}_1, \dot{q}_2)=  (1, 2, 0, 0)$, $(1, -1, 0, 0)$, 
$(-2, -1, 0, 0)$ and the origin. The first three are indicated by $C_1$, $C_2$, $C_3$. 
The eigenvalues near the origin are $-a-1 \pm \sqrt{a^2-a+1}$ corresponding with the frequencies 
$\omega_{2, 3}$, $\omega_{4, 5}$. For the normal modes in  system \eqref{IMsym6} we have one optical, and one acoustic. 
For $a=0.01$ the  eigenvalues (frequencies squared) are:
\[ \omega_1^2,\, \omega_2^2 =  {0.0149623,\, 2.00504} \] 
Denoting the critical points by $(q_1,  q_2)$ we find the eigenvalues: $(1, 2) \rightarrow \pm 2.4525\,i,\, \pm 0.1223$, 
$(1, -1) \rightarrow \pm 0.5477\,i,\, \pm 0.5477$, $(-2, -1) \rightarrow \pm 0.5758\,i,\, \pm 0.5211$. We conclude that for each of the nontrivial 
equilibria we have a centre manifold corresponding with 2 purely imaginary eigenvalues, one stable manifold and one unstable manifold. 
All nontrivial equilibria are unstable and may produce unbounded solutions.\\

\begin{figure}[ht]
\begin{center}
\includegraphics[width=7cm]{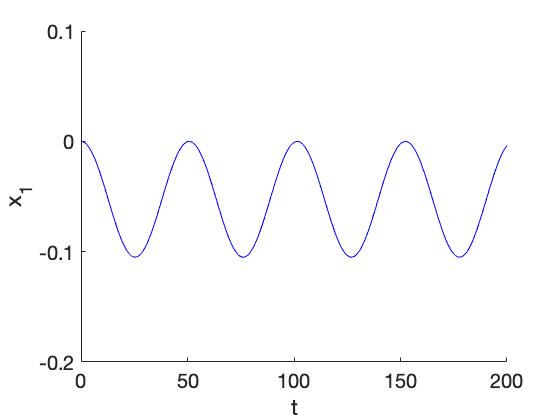}\, \includegraphics[width=7cm]{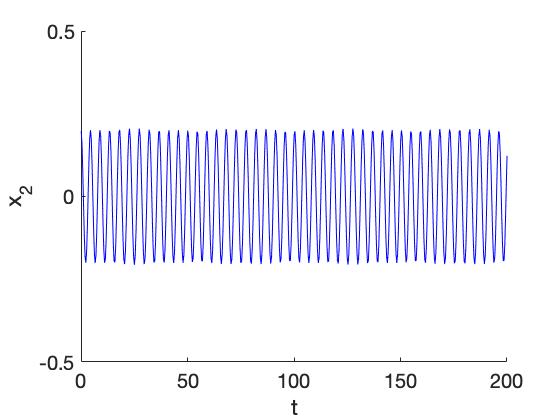}
\end{center}
 \caption{Solutions within the symmetric invariant manifold in the  case $p=3$ described by eq.~\eqref{eqsp=3a}, 
 $a=0.01$, $\alpha = 1$. 
 The coordinates are $x_1, v_1= \dot{x}_1$ (acoustic) and  $x_2, v_2= \dot{x}_2$ (optical)  of the  $\alpha$-chain 
 with 6 particles in 200 timesteps.  The initial values  are
 $x_1(0)= v_1(0)=0$, $x_2(0)= 0.2$, $v_2(0)=0$. We observe forcing of the acoustic mode. } 
\label{figp=3a}
\end{figure}   

Note that the full 6 dof Hamiltonian system contains system \eqref{IMsym6} supplemented 
by the (mirrored) modes $q_4, q_5$. The symmetry assumptions reduce the 6 dof system to two equivalent 
2 dof systems. 
Using the eigenvalues and eigenvectors for system \eqref{IMsym6} we can construct a 2 dof system in quasi-harmonic form 
with dynamics described by:
\begin{eqnarray} \label{eqsp=3a}  
\begin{cases}
 \ddot{x}_1 +  0.0149623x_1 & =   \alpha(0.0074805894 x_1^2 +  0.0004488580 x_1 x_2 -  0.0395000057 x_2^2), \\
\ddot{x}_2  +  2.00504x_2 & =   \alpha(0.0056956264 x_1^2 -  2.0048858034 x_1 x_2 -  0.0300748059 x_2^2) .
\end{cases}
\end{eqnarray} 
The 2 frequencies $\omega_1, \omega_2$ are not resonant; according to Lyapunov we can continue the linear normal modes 
in a neighbourhood of the origin. It is straightforward to obtain a convergent series approximation of the normal modes by 
the Poincar\'e-Lindstedt method. 

\subsection{Interactions in a neighbourhood of the origin} 

System \eqref{eqsp=3a} is, apart from the coefficients, symmetric in $x_1, x_2$ showing similar forcing of acoustic modes 
by optical ones and vice versa.  In the equation for $x_1$ the square $x_2^2$ will be forcing, in the equation for $x_2$ 
this will be $x_1^2$. 
Starting with initial zero values in the optical mode  the normal mode $x_1(t)$ 
shows interactions with the optical mode. Continuing the optical mode with zero acoustic initial values we 
find small interactions. \\
An asymptotic approximation scheme runs as follows. Near the origin we rescale the coordinates $x \mapsto \varepsilon x$ 
and divide the resulting equations by $\varepsilon$. We may replace $\varepsilon \alpha$ by $\varepsilon$ in system 
\eqref{eqsp=3a}. The lowest order normal form is easy to obtain, see for the procedure \cite{SVM}. The normal form has 
3 integrals, the Hamiltonian 
\begin{equation} \label{intp=3} 
H = \frac{1}{2}(\dot{x}_1^2 + \omega_1^2 x_1^2) + \frac{1}{2}(\dot{x}_2^2 + \omega_2^2 x_2^2) -\varepsilon 0.00249x_1^3 
+ \varepsilon 0.001x_2^3, 
\end{equation} 
and the 2 actions $\frac{1}{2}(\dot{x}_j^2 + \omega_j^2 x_j^2)$, $ j= 1, 2$. Quasiperiodic interaction will show at the next order 
of approximation. A more direct and easier approach runs as follows.\\
To approximate the quasi-periodic flow we use the lowest order $\varepsilon =0$ approximations: 
\begin{equation} \label{app0} 
x_1^0(t)= a_0 \cos \omega_1 t + b_0 \sin \omega_1 t, x_2^0(t)= c_0 \cos \omega_2 t + d_0 \sin \omega_2 t.
\end{equation} 
Introducing the zero order approximation \eqref{app0} in the right-hand sides of system \eqref{eqsp=3a} we obtain an 
approximation of first order; as the modes are not resonant we get no secular terms. 

Some numerical experiments for $\alpha = 1$ are shown in fig.~\ref{figp=3a} where we abbreviate in the caption $\dot{x}=v$. 
In fig.~\ref{figp=3a} we have $ a= 0.01$, $x_2(0)=0.2$ and the other initial values zero. We 
find interaction in the form of excitation of the acoustic mode $x_1$.\\ 
Consider now the first order approximation of $x_1(t)$ described above, If we start with the initial 
values used in fig.~\ref{figp=3a} we have for $t \geq 0$, $x_1^0(t)=0$,  $ x_2^0(t)= 0.2 \cos \omega_2 t $. 
The next step is solving an inhomogeneous linear equation to find the excitation of the acoustic mode $x_1(t)$ from: 
\[ \ddot{x}^1_1 +  0.0149623x^1_1 = -  0.0395000057 (0.2 \cos \omega_2 t)^2. \] 
Comparing $x^1_1(t)$ with the result in fig.~\ref{figp=3a} we find on 200 timesteps the same behaviour with error 0.02. 
Extending the Poincar\'e expansion procedure we can extend the accuracy and/or the time interval. Another 
possibility is to use 2nd order averaging. \\

In the case of 2 dof with resonant frequencies the phase-space dynamics near stable equilibrium is more interesting. 
In our case 
of widely separated frequencies we find quasiperiodic solutions corresponding with families of tori but as prominent 
feature significant interaction of modes. 

\subsection{Interactions near unstable equilibria} 

\begin{figure}[ht]
\begin{center}
\includegraphics[width=7cm]{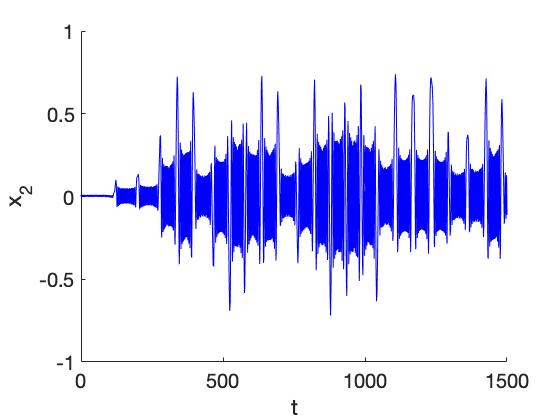}\, \includegraphics[width=7cm]{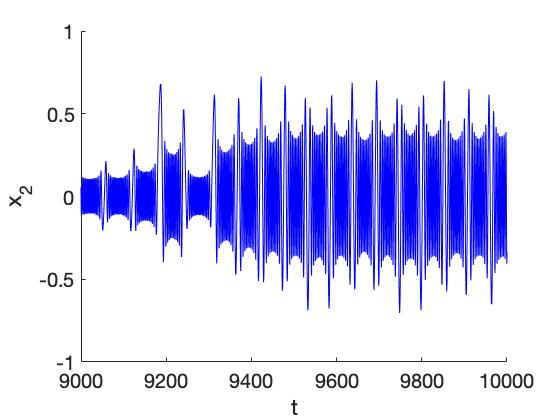}
\end{center}
 \caption{Solutions in the  case $p=3$ described by eq.~\eqref{eqsp=3a}, 
 $a=0.01$, $\alpha = 1$. 
 The coordinates are $x_1$, $v_1= \dot{x}_1$ (acoustic) and  $x_2$, $v_2= \dot{x}_2$ (optical)  of the  $\alpha$-chain 
 with 6 particles; left $x_2(t)$ in 1500 timesteps, right the time series from 9k till 10k.  The initial values  are
 $x_1(0)= 2$, $v_1(0)=0$, $x_2(0)= 0.01$, $v_2(0)=0$ near unstable equilibrium. 
We observe forcing of the optical mode $x_2$.   } 
\label{figp=3b}
\end{figure}   

The general dynamics merits closer attention, we perform calculations for larger initial values, take $\alpha =1$.
As an example we consider the dynamics near equilibrium $C_1$ with coordinate values $(1, 2, 0, 0)$ in system \eqref{IMsym6}. 
The equilibrium is unstable and although we can approximate the centre manifold with high accuracy, we can not be certain 
to start exactly in this submanifold. 
Using system \eqref{eqsp=3a} we start with initial conditions $(2, 0, 0, 0)$. Inverting our transformations we find that 
for system \eqref{IMsym6} this implies $q_1(0)= 1. 00754$, $q_2(0)=2$ so we start close to equilibrium $C_1$. The 
numerical solutions increase rapidly to very large numbers, suggesting unbounded behaviour.\\
Next we choose  initial conditions $(2, 0.01, 0, 0)$ corresponding in  system \eqref{IMsym6} with $q_1(0)= 0. 98768$, $q_2(0)=2.0001$, 
again close to $C_1$. The numerics shown in fig.~\ref{figp=3b} suggest that the initial values are 
located in a stable centre manifold; the time series $x_1(t)$ (not shown) shows oscillations between $-1$ and 2. We have  
strong excitation of the optical mode.  The stable and unstable manifolds of the 4-dimensional equilibrium are outside the 
center manifold, but can still be near the orbits. 
The light-grey vertical segments show relatively fast motion away from these 2-dimensional saddle structures.

\section{Periodic FPU $\alpha$-chain with 10  particles} \label{sec4} 

As the next case we consider  $N=10$. With $p=5$ the symmetry assumptions \eqref{IM2ps} produce invariant 
manifolds with dynamics described by 2 equivalent 4 dof systems. We have for the 1st system: 
\begin{eqnarray} \label{IM10} 
\begin{cases}
 \ddot{q}_1 +2 q_1  -q_2 & = q_2^2 -2q_1q_2,\\
m \ddot{q}_2 +2 q_2 - q_{1} -q_3 & = q_3^2 -2q_2q_3+2q_1q_2 -q_{1}^2,\\
 \ddot{q}_3 +2 q_3 - q_{2} -q_4 & = q_4^2 -2q_3q_4 +2q_2q_3 -q_{2}^2,\\ 
m \ddot{q}_4 +2 q_4 - q_{3}  & =   2q_3q_4-q_{3}^2.  
 \end{cases}
\end{eqnarray} 
The first two equations are identical to the first two of the general case in system \eqref{IM2p}, but the linear transformation to 
quasi-harmonic equations will be different in the general case. The eigenvalues of the linear part on the left-hand side are:
\begin{equation} \label{eigN=10} 
a+1 \pm \sqrt{1-\frac{1}{2}(1+ \surd{5})a +a^2},\,\ a+1 \pm \sqrt{1-\frac{1}{2}(1- \surd{5})a +a^2},
\end{equation} 

\begin{figure}[ht]
\begin{center}
\includegraphics[width=7cm]{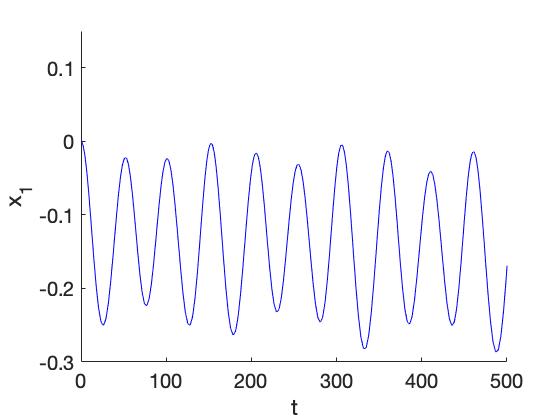}\,\includegraphics[width=7cm]{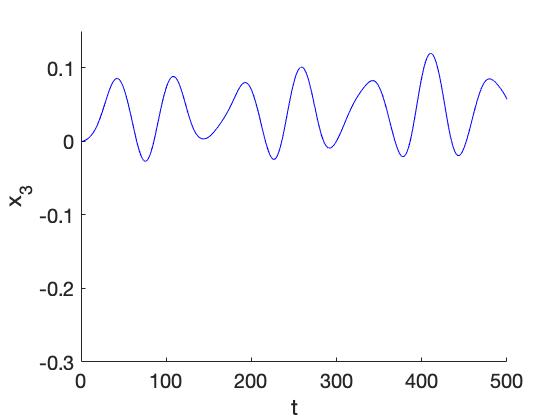} \\
\end{center}
 \caption{Solutions of system \eqref{eqsp=5} (case $p=5$) producing 2 invariant manifolds with 4 dof and equivalent dynamics, $a=0.01$. 
 Initial values $x_1(0), x_3(0)=0$ (acoustic )  and $x_2(0), x_4(0)=0.15$ (optical)   in 500 timesteps.  The initial   
 velocities are zero. 
We observe forcing of the acoustic modes and recurrence.} 
\label{figp=5}
\end{figure}  

\subsection{Equilibria} 
Equilibria correspond with critical points of the vector field describing system \eqref{IM10}. To find the equilibria 
we have to solve 4 quadratic algebraic equations. Apart from the trivial solution (the origin of phasespace) there are at most 
15 different solutions. With velocities zero we find easily: 
$(q_1, q_2, q_3, q_4) = (2, -1, 1, -2)$,  $(2, -1, 1, 3)$, $(3, -1, 1, -1)$, $(3, -1, 1, 3)$; more equilibria can be obtained  
by considering the cases $2q_1=q_2$ and $2q_3=q_4$. 

\subsection{Dynamics and interaction near the origin}
The origin corresponds with stable equilibrium. After transformation to quasi-harmonic form we consider the 
dynamics described by the system: 

\begin{eqnarray} \label{eqsp=5} 
\begin{cases}
\ddot{x}_1  + 0.0180728 x_1 & =  -0.04 x_4^2-0.13 x_2 x_4 +0.018 x_1 x_3+0.00065x_2x_3+0.009 x_3^2+  \\
& \hspace{1cm} 0.0004 x_1 x_4+0.00025 x_3 x_4,  \\

\ddot{x}_2 + 2.00193 x_2 & =  + 0.0063 x_1x_3-0.0019 x_3^2 - 2.0018 x_2 x_3-1.251 x_1 x_4-0.045 x_2 x_4+ \\
& \hspace{1cm} 1.251 x_3 x_4+0.0086 x_4^2, \\
 
\ddot{x}_3 +  0.00686474 x_3 & = -0.039x_2^2+0.049x_2x_4+ 0.00025 x_1x_2+0.0069 x_1x_3 -0.00015 x_2x_3 + \\
& \hspace{1cm}  0.00009 x_1x_4 + 0.0034x_1^2,\\  
 
\ddot{x}_4 + 2.01314 x_4  & = 0.0051 x_1^2 +0.0062 x_1 x_3  - 3.22 x_1 x_2 -0.05798 x_2^2 +3.22 x_2 x_3-\\ 
& \hspace{1cm}  2.013 x_1 x_4 +0.0443 x_2 x_4.
 \end{cases}
\end{eqnarray} 
The modes $x_1$, $x_3$ belong to the acoustic group, $x_2$, $x_4$ are optical. The largest coefficients of the quadratic terms 
are found in the optical group $x_2$, $x_4$. This will also influence the dynamics. \\
A difference with the case $p=3$ is that 
within the acoustic and optical group we have a detuned $1:1$ resonance. In \cite{FV79} the $1:1$ resonance is studied 
with results on in-phase and out-of-phase periodic solutions; in \cite{E05} (sections 1.3 and 3) a systematic study was set up, 
for instance regarding the $1:1:1:1$ resonance.  
We note that the optical modes in system 
\eqref{eqsp=5} are closer to exact resonance than the acoustic ones.   
For each equation we have 2 forcing terms, for $x_1$ they are $x_4^2$ and $x_2x_4$; for $x_2$ we have $x_3^2$ and $x_1x_3$. 
Terms like $x_2x_3$ or $x_1x_4$ have initially less influence as can be shown by averaging. \\
In fig.~\ref{figp=5} we start with nonzero optical modes $x_2$, $x_4$ leading to significant $x_1(t)$ oscillations.\\
\ \\
{\bf Invariant manifolds and normal modes}\\ 

\begin{figure}[ht]
\begin{center}
\includegraphics[width=7cm]{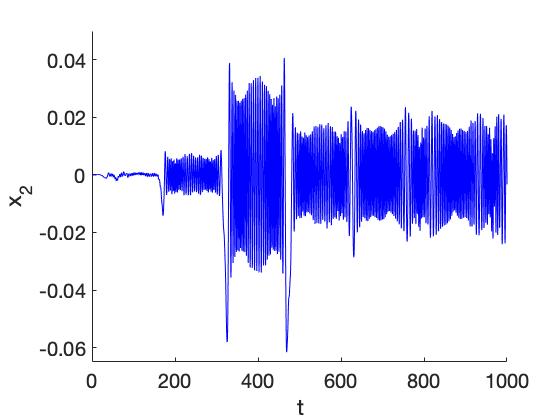}\,\includegraphics[width=7cm]{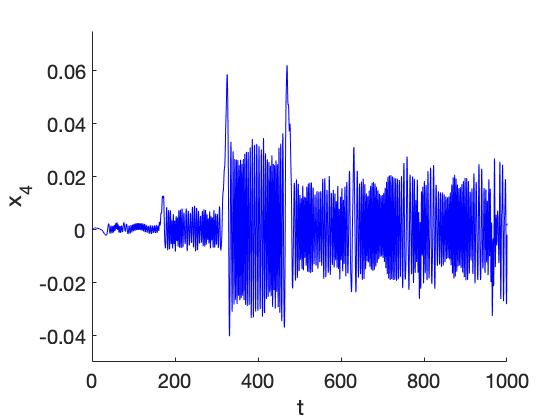} 
\end{center}
 \caption{Solutions in the case $p=5$ starting near the unstable $x_1$ normal mode in system \eqref{eqsp=5}. 
 The initial conditions are $x_1(0)=0.25$, $\dot{x}_1(0)=0.065$, $x_3(0)=0.01$, $\dot{x}_3=0.01$.
 The initial $x_2$, $x_4$ positions and velocities are zero. We observe excitation of the optical modes $x_2$, $x_4$.} 
\label{tenpartE1E4}
\end{figure}  

As the two $1:1$ resonances are detuned we can apply Lyapunov continuation, the normal modes exist. 
However, all the normal modes are unstable because of the high-low frequency interaction. \\
System \eqref{eqsp=5} contains 2 invariant manifolds described by the systems: 
\begin{equation} \label{subman14p=5}
\ddot{x}_1  + 0.01807 x_1  =   -0.04 x_4^2 + 0.0004 x_1 x_4,\  \ddot{x}_4 + 2.01314 x_4 = -0.0051 x_1^2 -  2.013 x_1 x_4,
\end{equation} 
\begin{equation}   \label{subman23p=5}
\ddot{x}_2 + 2.00193 x_2  =   -0.0019 x_3^2 - 2.0018 x_2 x_3,\ \ddot{x}_3 +  0.00686474 x_3 = -0.039x_2^2 -0.00015 x_2x_3.
\end{equation} 
Within these invariant manifolds the normal modes are harmonic and unstable but on long intervals of time.  Because of 
the recurrence theorem for measure-preserving maps, in particular Hamiltonian systems on a bounded energy manifold, 
high frequency modes that excite low frequency modes imply that the reverse  will also happen, except that the energy 
stored in low frequency modes is much less that in high frequency modes. We demonstrate this for the submanifold consisting 
of the $(x_1, x_4)$ modes of system \eqref{subman14p=5}. In fig.~\ref{tenpartE1E4} we show excitation of mode 4 when
starting close to the unstable $x_1$ normal mode.

\section{Periodic FPU $\alpha$-chain with 18  particles ($p$ not prime)} \label{sec5} 
The motivation to study briefly the case $p=9$ is to compare with the results in \cite{FVa} where after the case $N=4$ 
the case $N=8$ keeps the interaction of the chain with 4 particles but also shows new phenomena, in particular new  
invariant manifolds. In the case of 18 particles we expect from the embedding theorem a submanifold 
corresponding with 6 particles ($p=3$).\\ 
We will use now the theory of symmetric invariant manifolds from section \ref{sec2}. 
Transforming the 18 equations of motion induced by Hamiltonian (\ref{Hfpu}) to quasi-harmonic
form by $q \mapsto x$ we find 2 symmetric invariant manifolds with 8 2nd order equations. 
The 8 eigenvalues $\lambda_j$ are:

\begin{equation}\label{eigsp=9}
0.019391,\, 2.00061,\, 0.0149623,\, 2.00504,\, 0.00821511,\, 2.01178,\, 
0.00231905,\, 2.01768. 
\end{equation} 

The acoustic group corresponds with indices j = 1, 3, 5, 7, the optical group with even
indices. Note that from the embedding theorem we expect submanifolds corresponding with 6 and 12 particles, 
but now we are working in the subsystem \eqref{IM2p} corresponding with $p-1=8$ particles.  

\subsection{The reduced system}
The system in $x$-coordinates becomes: 

\begin{eqnarray} \label{eqsp=9} 
\begin{cases}
\ddot{x}_1  + \lambda_1 x_1 =  0.00970 {x_5^2}+0.0194 {x_3}
   {x_5}+0.000490 {x_4}
   {x_5}+0.000320 {x_6} {x_5}  -0.0401 {x_6^2} \\
 \hspace{1cm}  +0.000434
   {x_3} {x_6}
-0.0906 {x_4}
   {x_6}
+0.0194 {x_1}{x_7}   +0.000746 {x_2}
   {x_7}+0.0194 {x_3}
   {x_7}\\
 \hspace{1cm}+0.000259 {x_4}
   {x_7}+0.000404 {x_1} {x_8} -0.352 {x_2}
   {x_8}+0.000354 {x_3}
   {x_8}-0.140 {x_4}
    {x_8},  \\

\ddot{x}_2 +\lambda_2 x_2  =  -0.00140 {x_5^2}+0.00380
   {x_3} {x_5}-0.794 {x_4}
   {x_5}+0.703 {x_6}
   {x_5}\\
 \hspace{1cm}
+0.00579 {x_6^2}-0.703
   {x_3} {x_6}-0.0177 {x_4}
   {x_6}+0.00229 {x_1} {x_7}\\
 \hspace{1cm}
-2.00
   {x_2} {x_7}-0.00201 {x_3}
   {x_7}
+0.794 {x_4} {x_7}-1.08
   {x_1} {x_8}\\
 \hspace{1cm}-0.0417 {x_2}
   {x_8}+1.08 {x_3} {x_8}+0.0145
   {x_4} {x_8}, \\
 
\ddot{x}_3 +  \lambda_3 x_3 = 0.00748
   {x_3^2}+0.000449 {x_4}
   {x_3}-0.0395 {x_4^2}+0.0150
   {x_1} {x_5}\\
 \hspace{1cm}+0.000953 {x_2}
   {x_5}+0.000335 {x_1}
   {x_6}-0.176 {x_2} {x_6}+0.0150
   {x_1} {x_7}\\
 \hspace{1cm}-0.000504 {x_2}
   {x_7}+0.0150 {x_5}
   {x_7}+0.000114 {x_6}
   {x_7}+0.000273 {x_1}
   {x_8}\\
 \hspace{1cm}+0.272 {x_2}
   {x_8}+0.000176 {x_5}
   {x_8}-0.0954 {x_6}
   {x_8},\\  
 
\ddot{x}_4 + \lambda_4 x_4  = 0.00570 {x_3^2}-2.00
   {x_4} {x_3}-0.0301
   {x_4^2}+0.00959 {x_1}
   {x_5}\\
 \hspace{1cm}-5.05 {x_2} {x_5}-1.78
   {x_1} {x_6}-0.113 {x_2}
   {x_6}+0.00507 {x_1} {x_7}\\
 \hspace{1cm}+5.05
   {x_2} {x_7}-0.00327 {x_5}
   {x_7}+1.78 {x_6} {x_7}-2.74
   {x_1} {x_8}\\
 \hspace{1cm}+0.0921 {x_2}
   {x_8}+2.74 {x_5} {x_8}+0.0209
   {x_6} {x_8}.\\

\ddot{x}_5 + \lambda_5 x_5  = 0.00411
   {x_7^2}+0.00822 {x_3}
   {x_7}-0.0000708 {x_4}
   {x_7}+0.0000377 {x_8}
   {x_7}\\
 \hspace{1cm}-0.0404 {x_8^2}+0.00821
   {x_1} {x_3}+0.000523 {x_2}
   {x_3}+0.000208 {x_1}
   {x_4}\\
 \hspace{1cm}-0.109 {x_2}
   {x_4}+0.00821 {x_1}
   {x_5}-0.000386 {x_2}
   {x_5}+0.000135 {x_1}
   {x_6}\\
 \hspace{1cm}+0.0968 {x_2}
   {x_6}+0.0000967 {x_3}
   {x_8}+0.0592 {x_4}
   {x_8}.\\

\ddot{x}_6 + \lambda_6 x_6  = -0.000724 {x_7^2}+0.00371
   {x_3} {x_7}+2.27 {x_4}
   {x_7}+3.10 {x_8} {x_7}\\
 \hspace{1cm}+0.00711
   {x_8^2}+0.0109 {x_1}
   {x_3}-5.73 {x_2} {x_3}-2.27
   {x_1} {x_4}\\
 \hspace{1cm}-0.145 {x_2}
   {x_4}+0.00801 {x_1} {x_5}+5.73
   {x_2} {x_5}-2.01 {x_1}
   {x_6}\\
 \hspace{1cm}+0.0944 {x_2} {x_6}-3.10
   {x_3} {x_8}+0.0267 {x_4}
   {x_8}.\\

\ddot{x}_7 + \lambda_7 x_7  = 0.00116 {x_1^2}+0.0000892
   {x_2} {x_1}+0.00232 {x_3}
   {x_1}+0.0000310 {x_4}
   {x_1}\\
 \hspace{1cm}-0.0389 {x_2^2}-0.0000781
   {x_2} {x_3}+0.0309 {x_2}
   {x_4}+0.00232 {x_3}
   {x_5}\\
 \hspace{1cm}-0.0000200 {x_4}
   {x_5}+0.0000177 {x_3}
   {x_6}+0.0108 {x_4}
   {x_6}+0.00232 {x_5}
   {x_7}\\
 \hspace{1cm}-\left(6.91\times 10^{-6}\right)
   {x_6} {x_7}+0.0000106 {x_5}
   {x_8}+0.0148 {x_6}
   {x_8}.\\

\ddot{x}_8 + \lambda_8 x_8  =0.00214 {x_1^2}-3.73
   {x_2} {x_1}+0.00374 {x_3}
   {x_1}-1.48 {x_4} {x_1}\\
 \hspace{1cm}-0.0717
   {x_2^2}+3.73 {x_2}
   {x_3}+0.0498 {x_2}
   {x_4}+0.00242 {x_3} {x_5}\\
 \hspace{1cm}+1.48
   {x_4} {x_5}-1.31 {x_3}
   {x_6}+0.0113 {x_4}
   {x_6}+0.000942 {x_5}
   {x_7}\\
 \hspace{1cm}+1.31 {x_6} {x_7}-2.02
   {x_5} {x_8}+0.00601 {x_6}x_8.
 \end{cases}
\end{eqnarray} 

A first consequence of the forcing is that the acoustic and optical groups do not exist as separate invariant manifolds. To 
be explicit: 
each equation in system \eqref{eqsp=9} contains one quadratic forcing term of the other group and 3 quadratic mixed 
forcing terms of the other group, for instance for the 
first equation $x_6^2$ and $x_4x_6$, $x_2x_8$, $x_4x_8$. The largest coefficients of the quadratic terms are found in the 
optical group for the
$x_6$, $x_8$ equations. Analogous to the case of $p=3$ we have detuned $(1:1:1:1)$ resonances that will affect the dynamics. \\

{\bf The $(x_3, x_4)$ invariant manifold.}\\ 
There is an invariant submanifold given by $x_1=x_2=x_5=x_6=x_7=x_8=0$. The system is: 
\begin{eqnarray} \label{eqsp=9sub} 
\begin{cases}

\ddot{x}_3 +  0.01496 x_3 & = 0.00748 {x_3^2}+0.000449 {x_4} {x_3}-0.0395 {x_4^2},\\ 
 
\ddot{x}_4 + 2.00504 x_4  & = 0.00570 {x_3^2}-2.00 {x_4} {x_3}-0.0301 {x_4^2}.\\ 
  
 \end{cases}
\end{eqnarray} 
The system \eqref{eqsp=9sub} is identical (modulo numeric abbreviations) to system \eqref{eqsp=3a} for $p=3$. 
The presence of system \eqref{eqsp=9sub} as a submanifold in the case $p=9$ is predicted by the embedding theorem 
of \cite{BValt}, see also \cite{FVa}. We found only 1 invariant manifold for system  \eqref{eqsp=9}.\\
{\bf Normal modes}\\
As the $1:1: 1:1$ resonances are all detuned we can apply Lyapunov continuation, the normal modes exist. Even when 
using normal forms, applying detuned resonance, we can use the Weinstein \cite{W73} theorem for periodic normal mode 
solutions. However, all the normal modes are unstable because of the high-low frequency interaction. From system 
\eqref{eqsp=9} we can list the quadratic forcing of each mode between brackets; we have 
\[x_1 (x_6),\, x_2 (x_5),\, 
x_3 (x_4),\, x_4 (x_3),\, x_5 (x_8),\, x_6 (x_7),\, x_7 (x_2),\, x_8 (x_1). \] 
So acoustic $x_1$ excites optical $x_6$, $x_6$ excites acoustic $x_7$ etc. This circularity 
of excitations persists for more than 18 dof, it will be an important aspect of further study.

\section{Alternating FPU-chains with many particles} \label{sec6} 
Our {\sc Mathematica} programme described below produces by suitable linear transformations 
quasi-harmonic systems of equations where the normal modes can be identified. 
As in our studies of the cases $p=3, 5$, but considering now all cases $p$ prime, $p \leq 47$ we find out again  
 whether for instance quadratic optical terms $x_i^2$ arise as forcing terms in the quasi-harmonic form of the acoustic part 
 of the system. This turns out to be the 
 case for all prime numbers $p =3, 5, 7, \ldots, 47$. 
 We conclude that we have  interaction between the acoustic and optical group for the alternating FPU-chain with 
 large mass up to 104 particles and all multiples of these cases. 
  
 \subsection{The notebook plotinteractionp.nb} 
 
We developed and used a Mathematica notebook  {\bf plotinteractionp.nb}, see for details \cite{KVN20}. 
Given an odd $p$ it first sets up the system in $q$-variables (\ref{IM2p}) in symbolic form. As described in section \ref{sec2} 
it is of the form $\ddot{q} + B q = \alpha N(q)$ with $B$ a $(p-1) \times (p-1)$ matrix and $N(q)$ a quadratic vector function. 
We need to get to quasiharmonic
form. We know the eigenvalues in symbolic form of the matrix $B$, but the notebook finds them anyway, also in symbolic form.
Finding the corresponding matrix of
eigenvectors takes longer. As the matrix that is found becomes singular for $a=0$, we rescale the offending rows by a factor $1/a$.
Then we substitute high precision numbers before inverting the matrix. We are ready to produce the system in quasiharmonic
form in $x$-variables and can start plotting solutions. 

\subsection{Forcing squares on the right-hand sides} \label{forcing}
When checking where the squares occur on the right-hand sides of a system in $x$-coordinates like system \eqref{eqsp=5}, a remarkable pattern arises. To study this pattern we considered in detail 
chains with $N=2p$ particles for odd $p$,  $p\leq47$.  
Recall that the modes come in pairs consisting of an acoustic mode with variable 
$x_{2j-1}$ and an optical mode 
with variable $x_{2j}$. The  modes $2j-1$ and $2j$
belong together in the sense that $\lambda_{2j-1}$, $\lambda_{2j}$ are the eigenvalues of a matrix
$$\left(
\begin{array}{cc}
 2 a & 2 a \cos \left(\frac{\pi  j}{p}\right) \\
 2 \cos \left(\frac{\pi  j}{p}\right) & 2 \\
\end{array}
\right).$$
Note that $\lambda_{2j} +\lambda_{2j-1}=2+2a$.
Compare with \eqref{eqsp=3a}.
Now if one makes {\sc Mathematica} look where $x_{2i-1}^2$, $x_{2i}^2$ occur, then one finds that there is exactly one $j$ so that
they both occur in the right-hand side of the equations for $\ddot{x}_{2j-1} +\lambda_{2j-1} x_{2j-1}$ and $\ddot{x}_{2j} +\lambda_{2j} x_{2j}$. And $x_{2i-1}^2$, $x_{2i}^2$ occur nowhere else. Explicit examples are systems \eqref{eqsp=5} and \eqref{eqsp=9}. 

Let us write $j=\rho(i)$, so that $\rho$ is
the permutation of $\{1,\cdots,(p-1)/2\}$ that associates $j$ to $i$.
We have {\sc Mathematica} determine $\rho$ for all odd $p$ up to 47.  
 As pointed out to us by Stienstra \cite{S21} one has 
\begin{equation}\label{Jan}
\rho(i)=\min(2i,p-2i)
\end{equation} in all these cases.
If $p$ is prime, $5\leq p\leq47$, then $\rho$ has no fixed points.
On the other hand, when $p$ is divisible by three, we may have $i=j=\rho(i)$ as in system \eqref{eqsp=9sub}. But notice that in \eqref{eigsp=9}
the eigenvalues are ordered differently, making the formula $\rho(i)=\min(2i,p-2i)$ less apparent.
{\sc Mathematica} computes the cycle decomposition of $\rho$.
There is a connection between the cycle structure of $\rho$ and the existence of  certain invariant submanifolds of phase space. 
Indeed if $Y\subset\{1,\cdots,p-1\}$ is a nonempty proper subset of indices such that
$$V(Y)=\{ x_i=\dot x_i=0 \mbox{ for } i\in Y\}$$ is an invariant submanifold of phase space, then $2i-1\notin Y$ implies $2\rho(i)-1\notin Y$.
So $\{ i \mid 2i-1\in Y\}$ is a union of cycles  of $\rho$. (We consider fixed points as cycles too.) 
Observe that $\{ i \mid 2i-1\in Y\}=\{ i \mid 2i\in Y\}$. 
Not every union $Y$ of cycles will give an invariant manifold, as
one must also take mixed terms into account.
But the cycle decomposition of $\rho$ greatly simplifies the search for  $Y$ with $V(Y)$ invariant. 

One may also use plotting to search experimentally for suitable subsets $Y$. 
Simply take initial values such that initially a single mode is active and then look at the plots of all modes.
Clearly it suffices to try one initially mode for each cycle of $\rho$. \\
For instance, when $p=17$ the permutation $\rho$ has two cycles and
none of them yields an invariant proper submanifold. In fact we did not
find any invariant manifold of type $V(Y)$ for primes $3\leq p\leq47$. We present a few results without technical details.\\

{\bf Multiples of $p=3$}\\ 
If $p=9$ (18 particles), see also section \ref{sec5}, $\rho$ has cycles of length 1, 3 and we recover subsystem \eqref{eqsp=3a} with 2 modes.\\
If $p=15$ (30 particles) we recover the submanifolds with 2 modes and 4 modes  respectively; in this case 5 is also a divisor.\\
 For $p=27$ (54 particles) we 
 find one $V(Y)$ with 8 modes, containing another one with 2 modes. \\
For $p=45$ we find four different invariant $V(Y)$. They have 2,  4, 8  or 14 modes, corresponding with the respective divisors 3, 5, 9, 15 of 45.
There are four containments:
The one with 14 modes contains both the one with 2 modes and the one with 4 modes. The one with 2 modes is contained in both the one with 8 modes and the one with 14 modes. We recognise the divisor relations between 3, 5, 9, 15. \\

{\bf Multiples of $p= 5$}\\
The case $p=5$ (10 particles) is discussed in section \ref{sec4}.\\
The case $p=15$ has also 3 as a divisor, see above.\\
If $p=25$ (50 particles) we find a $V(Y)$ with 4 modes. \\

More generally, if $q$ is a proper divisor of an odd $p$, $p\leq47$,
we find a subsystem of $(q-1)$ modes inside the system in $x$-coordinates
that has $(p-1)$ modes. More specifically, the $(q-1)$ modes that remain 
are the ones involving $x_i$
whose $\lambda_i$ is an eigenvalue of a matrix
$$\left(
\begin{array}{cc}
 2 a & 2 a \cos \left(\frac{\pi  j}{q}\right) \\
 2 \cos \left(\frac{\pi  j}{q}\right) & 2 \\
\end{array}
\right)$$
with $1\leq j\leq (q-1)/2$.
So that is the form that the embedding theorem \ref{embed} now takes.\\

For odd $p$, $p\leq 47$, we find no $V(Y)$ that are not explained by the embedding theorem.

\subsection{Quadratic terms and normal forms} 
\begin{figure}[ht]
\begin{center}
\includegraphics[width=7cm]{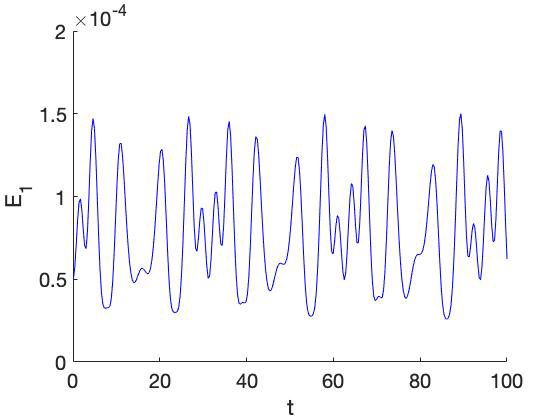}\,\includegraphics[width=7cm]{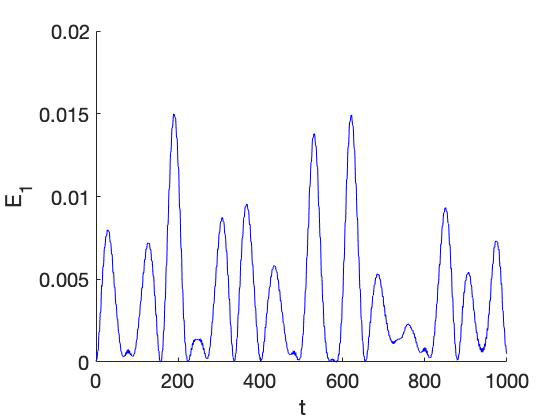} 
\end{center}
 \caption{The action $E_1(t)$ for 2 cartoon problems. Left the behaviour in system \eqref{car1} where 
 $\omega_1=0.1$, $\omega_2=1$, $x(0)=0.1$, $y(0)=0.1$, $E_1(0)= 0.00005$ and initial velocities zero. 
Right the behaviour in system \eqref{car2} where 
 $\omega_1=0.1$, $\omega_2=1$, $\omega_3=1.05$, $x(0)=0.1$, $y(0)=0.3$, $z(0)=0.3$, $E_1(0)= 0.00005$ and initial 
 velocities zero.  } 
\label{cartoon}
\end{figure}  

Near stable equilibrium it is natural to apply normal form analysis to the systems. Apart from forcing squares  we have 
mixed terms of the acoustic 
and optical groups but these terms are non-resonant; they will appear in the normal form at very high order and have little 
influence. Forcing quadratic terms consisting of only mixed acoustic or only mixed optical contributions will have 
forcing influence as they are in detuned resonance. 
For instance in the case of system  \eqref{eqsp=9} 
a calculation shows that we are still left with 8 quadratic terms 
on the right-hand side of each equation. Analysing the normal form poses a formidable problem.\\
Instead of normalising the systems obtained in the preceding sections we demonstrate the influence of mixed 
quadratic terms by integrating 2 typical Hamiltonian cartoon problems. In both cases we have 
coupled oscillators with widely separated frequencies; see \cite{TV} 
for an introduction and references. The standard procedure would be to introduce slowly-varying variables, for instance 
amplitude-phase or action-angle variables, see \cite{SVM}, to compute by averaging-normalisation an approximating 
normal form system. We leave out these technical details and will discuss the numerical results shown in fig.~\ref{cartoon}. 
The first cartoon is the system: 
\begin{eqnarray} \label{car1} 
\begin{cases} 
\ddot{x} + \omega_1^2 x = xy,\\
\ddot{y} + \omega_2^2 y = \frac{1}{2}x^2.
\end{cases} 
\end{eqnarray} 
We introduce the action $E_1= 0.5 (\dot{x}^2 + \omega_1^2 x^2)$. In fig.~\ref{cartoon} (left) we have $0< \omega_1 \ll 
\omega_2$. As predicted by normalisation the action $E_1(t)$ shows only small variations around its initial value.  
In the 2nd cartoon we have detuned forcing: 
\begin{eqnarray} \label{car2} 
\begin{cases} 
\ddot{x} + \omega_1^2 x = 0.2 yz,\\
\ddot{y} + \omega_2^2 y = 0.2 xz+0.25(z^2 + 2yz),\\
\ddot{z} + \omega_3^2 z = 0.2 xy+0.25(y^2 + 2yz),
\end{cases} 
\end{eqnarray} 
where $\omega_1 \ll \omega_2$ and $\omega_2$ is close to $\omega_3$. The $x$ normal mode is harmonic.
The $y, z$ oscillators are in detuned resonance and are strongly forcing the $x$-oscillator, see fig.~\ref{cartoon} (right). 

\section{Conclusions}  
\begin{enumerate} 

\item We have demonstrated interaction between acoustic and optical modes for periodic FPU-chains with alternating large mass  up to 104 particles and their multiples.

\item It was shown in section 5 of \cite{BValt}  that if $N=8$ we find 3 invariant manifolds with 3 dof in the case of the periodic FPU $\alpha$-chain with alternating masses. One of them corresponds with the dynamics of an alternating FPU-chain 
with 4 particles. The case $N=18$ in our section \ref{sec5} shows a different structure of invariant manifolds.\\
Considering the case of $N=2p$ particles with $p$ prime or odd we have obtained some general insight in the existence and structure 
of the invariant manifolds in the systems described by system \eqref{IM2p}. 

\item An open problem was formulated in section \ref{sec6} regarding the presence of quadratic terms in the 
quasi-harmonic 
form of the FPU-chains considered here. The circularity of exciting quadratic terms noted in 
section \ref{sec5} plays probably an essential part.
A general solution to this problem might answer the interaction question for 
systems with an arbitrary number of even particles. 

\item An interesting open problem is the high-low frequency interaction problem for $N=2p$ with $p$ prime, $p \geq 53$. 
The form of the invariant manifolds described by system \eqref{IM2p} in section \ref{sec2} suggest that such interactions 
take place. 

\end{enumerate} 

{\bf Acknowledgement} A useful remark by dr.\ J. Stienstra helped us in the discussion of permutations of quadratic terms in subsection \ref{forcing}. J. Stienstra also checked the transformations to equations of motions 
independently. We owe ref.~\cite{B20} to one of the anonymous referees.


\begin{thebibliography}{99} 

\bibitem{ADR} I.V. Andrianov, V.V. Danishevskyy and G. Rogerson, {\em Vibrations of nonlinear elastic lattices: low- and 
high-frequency dynamic models, internal resonances and modes coupling}, Proc. R. Soc. A 476: 20190532, 
\href{https://doi.org/10.1098/rspa.2019.0532}{DOI: 10.1098/rspa.2019.0532}

\bibitem{B20} G. Benettin and A. Ponno, {\em Understanding the FPU state in FPU-like models}. Mathematics in 
Engineering 3, pp. 1-22, \href{https://doi.org/10.3934/mine.2021025}{DOI: 10.3934/mine.2021025} (2020) 

\bibitem{BS12}
T. Bountis and H. Skokos, {\em Complex Hamiltonian Dynamics}, Springer (2012). 

\bibitem{BValt} R.W. Bruggeman and F. Verhulst, {\em Near-integrability and recurrence in FPU 
chains with alternating masses}, J. Nonlinear Science 29, pp. 183-206, 
\href{https://doi.org/10.1007/s00332-018-9482-x}{DOI: 10.1007/s00332-018-9482-x} (2019).

 \bibitem{BVapp} R.W. Bruggeman and F. Verhulst, {\em FPU-chains 
 with large alternating masses (plus appendix)}, \href{http://arxiv.org/abs/2004.03876}{arXiv 2004.03876} (2020).
 
\bibitem{CRZ05} D.K. Campbell, P. Rosenau and G.M. Zaslavsky (eds.), {\em The Fermi-Pasta-Ulam Problem. 
The first 50 years.}
Chaos, Focus issue 15, \href{https://doi.org/10.1063/1.1889345}{DOI: 10.1063/1.1889345} (2005). 

\bibitem{CEB} H. Christodoulidi, Ch. Efthymiopoulos,  and T. Bountis,  [2010]
{\em Energy localization on $q$-tori, long-term stability, and the
interpretation of Fermi-Pasta-Ulam recurrences}, Physical Review E 81,6210.

\bibitem{C05} G.M. Czechin., D.S. Ryabov and K.G. Zhukov, 
{\em Stability of low-dimensional bushes of vibrational modes in the Fermi-Pasta-Ulam chains}, 
Physica D 203, pp. 121-166 (2005). 
 
  \bibitem{DLS18} A. Degasperis, S. Lombardo and M. Sommacal, {\em Integrability and linear stability of nonlinear 
waves)}, J. Nonlinear Science 28, pp. 1251-1291 (2018).

\bibitem{E05} K. Efstathiou, {\em Metamorphoses of Hamiltonian systems with symmetries}, Lecture Notes in 
Mathematics 1864, Springer (2005). 

\bibitem{G08} G. Gallavotti (ed.) {\em The Fermi-Pasta-Ulam Problem: a status report}, Lecture Notes in 
Physics, Springer (2008). 

\bibitem{GGMV} L. Galgani, A. Giorgilli, A. Martinoli and S. Vanzini, On the problem of energy partition for large 
systems of the Fermi-Pasta-Ulam type: analytical and numerical estimates, Physica D 59, pp. 334-348 (1992).  

\bibitem{PR} P. Poggi  and S. Ruffo  {\em Exact solutions in the FPU oscillator chain}, Physica D 103, 
pp. 251-271 (1997). 

\bibitem{SVM} J.A. Sanders, F. Verhulst and J. Murdock, {\em Averaging methods in nonlinear dynamical systems}, 
2nd ed., Appl. Math. Sciences 59, Springer (2007).  

\bibitem{S21} J. Stienstra, Department of Mathematics, University of Utrecht, The Netherlands, {\em Private 
Communication} (2021).

\bibitem{TV} J.M. Tuwankotta and F. Verhulst, {\em Hamiltonian Systems with Widely Separated Frequencies}, 
Nonlinearity 16, pp. 689-706 (2003).

\bibitem{KVN20}  Wilberd van der Kallen and Ferdinand Verhulst, {\em  Explorations for alternating FPU-chains with large mass} (Includes Notebook), 
\href{http://arxiv.org/abs/2011.04290}{arXiv 2011.04290} (2020). 

\bibitem{FV79} F. Verhulst, {\em Discrete symmetric dynamical systems at the main resonances with applications 
to axi-symmetric galaxies}, Phil. Trans. roy. Soc. London 290, pp. 435-465 (1979). 

\bibitem{FV18}  Ferdinand Verhulst, {\em  Interaction of lower and higher order Hamiltonian resonances}, 
Int. J. Bif. Chaos 28, \href{https://doi.org/10.1142/S0218127418500979}{DOI 10.1142/S0218127418500979} (2018).


\bibitem{FVa}  Ferdinand Verhulst, {\em  High-low frequency interaction in alternating FPU  $ \alpha$-chains},  
Int. J. Nonlinear Mechanics, \href{https://doi.org/10.1016/j.ijnonlinmec.2021.103686}{DOI: 10.1016/j.ijnonlinmec.2021.103686} (Febr.~2021).

\bibitem{W73} A. Weinstein, {\em Normal modes for nonlinear Hamiltonian systems}, Invent. Math. 20, pp. 47-57 (1973).

\end{thebibliography}
\end{document}